\documentclass[12pt]{amsart}
\usepackage{graphicx,amssymb}
\usepackage{mathrsfs}
\usepackage{amssymb,amsmath,amsthm,color}
\usepackage{graphicx,mcite}
\usepackage{hyperref}
\usepackage{url}
\usepackage{setspace}

\theoremstyle{definition}

\theoremstyle{remark}

\numberwithin{equation}{section}

\newcommand{\R}{\mathbb R}

\def\TagOnRight

\def\R {\mathbb{R}}

\newcommand{\be}{\begin{equation}}
\newcommand{\ee}{\end{equation}}
\newcommand{\bea}{\begin{eqnarray}}
\newcommand{\eea}{\end{eqnarray}}
\newcommand{\Bea}{\begin{eqnarray*}}
\newcommand{\Eea}{\end{eqnarray*}}
\newcommand{\bt}{\begin{Theorem}}
\newcommand{\et}{\end{Theorem}}

\newcommand{\bpr}{\begin{Proposition}}
\newcommand{\epr}{\end{Proposition}}

\newcommand{\bl}{\begin{Lemma}}
\newcommand{\el}{\end{Lemma}}
\newcommand{\bi}{\begin{itemize}}
\newcommand{\ei}{\end{itemize}}

\newtheorem{Definition}{Definition}[section]
\newtheorem{Theorem}[Definition]{Theorem}
\newtheorem{Lemma}[Definition]{Lemma}
\newtheorem{Proposition}[Definition]{Proposition}

\newtheorem{Remark}[Definition]{Remark}

\begin{document}
\baselineskip16pt

\title{Maximal averages associated with families of finite type surfaces}
\author{Ramesh Manna}
\address{School of Mathematics, Harish-Chandra Research Institute, Allahabad, India 211019.}
\address {\& Homi Bhabha National Institute, Training School Complex, Anushakti Nagar, Mumbai 400085, India.}
\email{rameshmanna@hri.res.in}
\date{\today}

\begin{abstract}

We study the boundedness problem for maximal operators $\mathcal{M}$ associated to averages along the families of hypersurfaces $S$ of finite type in $\mathbb{R}^n.$ In this paper, we prove that if $S$ is a finite type hypersurface which is of finite type $k$ at $x_0 \in \mathbb{R}^n$, then the associated maximal operator is bounded on $L^p(\mathbb{R}^n)$ for $p>k.$ We shall also consider a variable coefficient version of maximal theorem and obtain the same $L^p$- boundedness result for $p>k.$

We also discuss the consequence of this result. In particular, we verify a conjecture by E. M. Stein and its generalization by A. Iosevich and E. Sawyer on the connection between the decay rate of the Fourier transform of the surface measure on the hypersurface $S$ and the $L^p$- boundedness of the associated maximal operator $\mathcal{M}.$  

\end{abstract}
\maketitle

\section{Introduction} \label{section1}
Let $S$ be a smooth hypersurface in $\mathbb{R}^n$ and let $\rho\in C_c^{\infty}(S)$ be a smooth function with compact support. Given a function $f,$ continuous and compactly supported, we consider for each $x \in \mathbb{R}^n$ and $t>0,$ the averaging operator 
$$M_tf(x) \, := \, \int_S f(x-ty) \, \rho(y) \, d\sigma(y),$$
where $d\sigma$ denotes the normalized Lebesgue measure over the hypersurface $S$.
Define the corresponding maximal operator by 
\begin{eqnarray} \label{M1}
\mathcal{M}f(x) \, := \, \sup_{t>0} |M_tf(x)|.
\end{eqnarray}

It is not obvious that such averaging operators are well defined for $f$ in $L^p-$ spaces, since $S$ has measure zero in $\R^n.$ Nevertheless, a priori $L^p$- estimates are possible when $S$ has suitable curvature properties.
Therefore, a natural question, we ask is for what range of the exponents $p$ is  the following a priori inequality satisfied:
\begin{eqnarray} \label{1.2}
\|\mathcal{M}f\|_{L^p(\mathbb{R}^n)} \leq \, B_p \, \|f\|_{L^p(\mathbb{R}^n)}, ~ f \in \mathscr{S}.
\end{eqnarray}

The aim of this paper is to study the $L^p$-mapping property of maximal operator associated with finite type hypersurfaces. There is a vast literature on maximal and averaging operators over families of lower dimensional surfaces of $\mathbb{R}^n.$ It turns out that curvature condition plays a crucial role in the analysis of $\mathcal{M}$: roughly speaking, curved surfaces admit non-trivial maximal estimates; whereas flat surfaces do not. A fundamental and representative positive result in this direction is the $L^p,~p>2$ boundedness of Bourgain's circular maximal operator.

The study of such maximal operator over dilations of a fixed hypersurface $S \subset \mathbb{R}^n$ has its beginnings in the work of E. M. Stein on the spherical maximal operator (see, e.g., \cite{ES}, \cite{EW}).
Stein showed that when $ S=\mathcal{S}^{n-1},$ the unit $n-1$ dimensional sphere, the corresponding maximal operator is bounded on $L^p(\mathbb{R}^n)$ if and only if $p > \frac{n}{n-1}, ~~ n \geq 3.$ Stein's proof of the spherical maximal theorem exploits curvature via the decay of the Fourier transform of the surface measure of the sphere. Indeed, one obtains the same sharp $L^p$- bounds if the sphere is replaced by a piece of any hypersurface in $\mathbb{R}^n$ with everywhere non-vanishing Gaussian curvature (see, \cite{G}).

The $2$- dimensional version of the spherical maximal operator was proved by Bourgain (see, \cite{B}). Bourgain's proof of the circular maximal theorem relies more directly on the geometry involved. The relevant geometry information concerns intersections of pairs of thin annuli, (for more details, see \cite{B}). It is well known and easy to see that, the maximal operator $\mathcal{M},$ given by \eqref{M1} will not be bounded on any $L^p$ for $p < \infty,$ if we consider the maximal averages over the boundaries of cubes instead of spheres. More generally, maximal averages over any hypersurface containing a piece of a hyperplane not containing the origin can never be bounded on $L^p$ for finite $p$. On the other hand, if the hypersurface is of finite type, then Sogge et. al. \cite{SS} showed that there exists a $p_0 < \infty$ such that inequality \eqref{1.2} holds for $p > p_0$.

In this paper, we shall consider a situation when the curvature is allowed to vanish of finite order on a finite set of isolated points. In dimension two,  Iosevich \cite{IOSE} considered the finite-type curve $\mathbb{C}$, which is of finite type $k$ at $x_0 \in \mathbb{R}^2.$ He proved that, the  corresponding maximal operator $\mathcal{M}$ satisfies the following inequality: 
$$\|\mathcal{M}f\|_{L^p(\mathbb{R}^n)} \leq \, B_p \, \|f\|_{L^p(\mathbb{R}^n)},$$  for $p > k,$ and also shown that this result is sharp. 
Here, we shall extend this result in higher dimension. We shall also consider a variable coefficient version of maximal theorem and obtain the same $L^p$- boundedness result for $p>k,$ see section 2. The proof of our main result, Theorem \eqref{Theorem1.1} will strongly make use of the results   of Iosevich \cite{IOSE} and Sogge \cite{CDS1}.

\section{Some Preliminaries and Main result}
We shall need the following definition: 
\begin{Definition}
We say that a smooth hypersurface $S$ is of finite type if $S$ has a finite order of contact with any affine hyperplane.
\end{Definition}
We shall also need a more precise definition which would specify the order of vanishing at each point. We consider $S$ in a sufficiently small neighbourhood of a given point and write $S$ as the image of a smooth mapping $\phi:U\rightarrow \mathbb{R}^n,$ where $U$ is a neighbourhood of the origin in $\mathbb{R}^{n-1}.$ Now, we fix a point $x_0 \in U,$ and a unit vector $\eta \in \mathbb{R}^n.$
We assume that the function $$[\phi(x)-\phi(x_0)]. \eta$$ does not vanish to infinite order as $x \rightarrow x_0,$ that is, 
for $x_0 \in U$ and each unit vector $\eta \in \mathbb{R}^n,$ there is a multi-index $\alpha,$ with $|\alpha| \geq 1,$ so that $$\left. \partial_x^{\alpha}[\phi(x). \eta] \right\vert_{x=x_0} \neq 0.$$
\begin{Definition}
Let $S$ be defined as before. We called the hypersurface $S$ is of finite type of order $k$ at $x_0$ if we have the smallest $k$ so that for each unit vector $\eta,$ there exists an $\alpha$ with $|\alpha| \leq k$ for which 
$$\left. \partial_x^{\alpha}[\phi(x). \eta] \right\vert_{x=x_0} \neq 0.$$
\end{Definition}
\begin{Remark}
In fact, this condition is equivalent to the condition that at least one of the principal curvatures of $S$ does not vanish to infinite order at $x_0.$
\end{Remark}

Let us recall at this point a result by A. Greenleaf \cite{G}. He proved that if 
\begin{eqnarray} \label{2.1}
\widehat{\rho d\sigma}(\xi) \, = O(|\xi|^{-\beta}) \mbox{ as } |\xi| \rightarrow \infty
\end{eqnarray}
 and if $\beta > \frac{1}{2},$ then the maximal operator is bounded on $L^p$ whenever $p> 1+ \frac{1}{2 \beta}.$ 

For $\beta= \frac{1}{2},$ E. M. stein and later for the full range $\beta \leq \frac{1}{2},$  A. Iosevich and E. Sawyer (see, \cite{IS}) conjectured that if $S$ is a smooth, compact hypersurface in $\mathbb{R}^n$ such that $\widehat{\rho d\sigma}(\xi) \, = O(|\xi|^{-\beta})$ for some $0 < \beta \leq \frac{1}{2},$ then the maximal operator $\mathcal{M}$ is bounded on $L^p(\mathbb{R}^n)$ for every $p > \frac{1}{\beta},$ at least if we assume $\rho > 0.$

A partial confirmation of stein's conjecture has been given by C. D. Sogge (see, \cite{CDS1})  who showed that the $L^p$- boundedness of the corresponding maximal operator $\mathcal{M}$ holds with one principal curvature nonzero everywhere on the surface, for every $p> 2.$ Also, we know that if the surface has at least one non-vanishing principal curvature than the estimate \eqref{2.1} above holds for $\beta= \frac{1}{2}$ (see, Littman \cite{L}).

Now, if $S$ is a smooth hypersurface in $\mathbb{R}^n$ of finite type of order $k$ at $x_0$. Then the estimate \eqref{2.1} holds for $\beta= \frac{1}{k},$ where $k$ is the type of $S$ inside the support of $\rho,$ (see, \cite{EMS}). Therefore, our result, Theorem \eqref{thm2.1} will give the confirmation of the Stein-Iosevich Sawyer conjecture for the finite type smooth-hypersurface.
\begin{Remark}
Notice that, the decay rate of the oscillatory integrals (Fourier transform of surface carried measure on the surface $S\subset \mathbb{R}^3$ of finite type of order $k$ at $x_0$) in $\mathbb{R}^3,$ given by $\beta=\frac{1}{k}$ (see, \cite{EMS}) is not sharp. In fact, the sharp decay of the Fourier transforms of surface carried measure on the surface $S\subset \mathbb{R}^3$ of finite type can be defined in terms of Newton polyhedrons (for more details, see, \cite{VAN, MIK}). 
\end{Remark}

Next, we shall state our main result of this paper. 

\begin{Theorem} \label{Theorem1.1}
Let $S$ be a finite-type hypersurface which is of finite type of order $k$ at $x_0$ in $\mathbb{R}^n.$ Let $$M_tf(x) \, = \, \int_S f(x-ty) \, \rho(y) \, d\sigma(y),$$  where $d\sigma$ is the induced Lebesgue measure on $S$ and a smooth cut off function $\rho$ supported in a sufficiently small neighbourhood of $x_0.$
Let, $$\mathcal{M}f(x) \, = \, \sup_{t>0} |M_tf(x)|.$$
Then, the inequality 
\begin{eqnarray} \label{thm2.1}
\|\mathcal{M}f\|_{L^p(\mathbb{R}^n)} \leq \, B_p \, \|f\|_{L^p(\mathbb{R}^n)},
\end{eqnarray}
 holds for $p > k,$ with some constant $B_p,$ depending only on $p.$
\end{Theorem}

We shall also consider the variable coefficient version of maximal Theorem \eqref{Theorem1.1}.

\subsection{Variable coefficient maximal theorem}
In this section, we shall discuss the variable coefficient version of maximal operator \eqref{M1}. The averaging operator $M_t$ that we have considered so far is called translation invariant or the constant coefficient operator, because it averages a function over the translates and dilates of a fixed surface. We are now going to consider an operator which averages a function over a more general distribution of surfaces in $\mathbb{R}^n,$ and also a more general time dependence. Iosevich introduced this operator in the plane.  As before, we are going to define a maximal operator by taking the supremum over the time dependence.

Let us recall at this point some of the previous results in this direction. Greenleaf \cite{G} considered the hypersuface of nonzero Gaussian curvature, and allow the surfaces to vary smoothly from point to point and behave only asymptotically like dilations. He proved that, for each compact set $K \subset \mathbb{R}^n$, the corresponding maximal operator is $L^p$- bounded for $p > \frac{n}{n-1}, ~ n \geq 3.$ Later C. D. Sogge showed that the same is true for $n=2$ (see \cite{S}).

In order to formulate the conditions for $L^p$ boundedness of such operators, we introduce a certain curvature hypothesis on the geometry of the surfaces. To discuss the curvature condition, let us first consider a surface distribution in $\mathbb{R}^n$ such that through every point $x \in \mathbb{R}^n,$ we have a smooth surface. We assume that this surface distribution is a smooth $n-1$ dimensional submanifold for each fixed $t.$ We then locally express this surface  distribution as $$D_t \, = \, \{(x, y) \in \mathbb{R}^n \times \mathbb{R}^n: x_1- y_1=A(x, y',t)\},$$ for some real and smooth $A(x, y', t), ~ y' =(y_2, \dots, y_n).$

Let $S_{x,t} \, = \, \{y \in \mathbb{R}^n: (x,y) \in D_t\}$ and consider a family of operators 
\begin{eqnarray} \label{4.1}
M_t(f)(x) \, = \, \int_{S_{x,t}} f(x-y) \, d \sigma_{x,t}(y),
\end{eqnarray}

where $d\sigma_{x,t}$ denotes the smooth cut off function $\psi(x,y,t)$ times the Lebesgue measure on $S_{x,t}.$ Let $\Phi(x,y,t)=x_1-y_1-A(x, y',t).$ Now, if we use the same notation by replacing $\psi$ by $|\nabla \Phi| \, \psi$ and then, using the defining function $\Phi,$ we can rewrite $M_t(f)(x)$ as 
\begin{eqnarray}
M_tf(x)&=&\int\limits_{\mathbb{R}^n}\delta_0(\Phi)f(x-y) \, \psi(x,y,t) \, dy \nonumber \\ 
&=&\int\limits_{\mathbb{R}^n} \, \delta_0 \left(x_1-y_1-A(x, y',t)\right) \, \psi(x, y, t) \, f(y) \, dy \nonumber \\
&=&\frac{1}{2\pi}\int\limits_{\mathbb{R}^n}\int\limits_{-\infty}^{\infty} \, e^{i\tau (x_1-y_1-A(x, y',t))} \, \psi(x, y, t) \, f(y) \, d\tau \, dy.
\end{eqnarray}

After setting $M_t(f)(x) : =Mf(x,t)$ if we regard $t$ as fixed, then it is well known that $M_t: D'(\mathbb{R}^n) \rightarrow D'(\mathbb{R}^n)$ is a Fourier integral operator of order $-\frac{n-1}{2},$ (see, \cite{SSS}). Let $k_t(x,y)$ denote the kernel of this Fourier integral operator. Then, the wave front set of a distribution defined by $(x,y) \rightarrow k_t(x,y)$ is contained in a subset of the cotangent bundle of $\mathbb{R}^n \times \mathbb{R}^n$ with zero section removed (see, e.g., \cite{S}, \cite{D}). This subset is called a canonical relation.

Let $C_t$ denote the canonical relation for a fixed $t,$ and also let $X, Y$ be the support of $x \rightarrow k_t(x, y)$ and $y \rightarrow k_t(x, y)$ respectively. Then, $C_t$ can be viewed as a Lagrangian submanifold of the cotangent space of $X \times Y$ with the zero section removed, where this is endowed with the usual canonical symplectic form $dx \wedge d \xi - dy \wedge d \mu$. Here, $\xi$ and $\mu$ are dual variable of $x$ and $y$ respectively. Let $\pi_l: C_t \rightarrow T^{\star}(X) \setminus \{0\}$ and $\pi_r: C_t \rightarrow T^{\star}(Y) \setminus \{0\}$ denote the natural projections of $C_t$, where, $T^{\star}X$ and $T^{\star}Y$  of course denote the co-tangent bundle of $X$ and $Y.$ 

We say that $C_t$ is a local canonical graph if $$C_t \, = \, \left\{(x, \xi, y, \mu): (y, \mu) = \chi_t(x, \xi)\right\},$$ where $\chi_t$ is a symplectomorphism for each $t.$ 
Notice that this condition is equivalent to the condition that $\pi_l$ and $\pi_r$ are local diffeomorphisms. 
In this paper, we shall always assume that the canonical relations $C_t$ of the operators are locally the graph of a canonical transformation. This condition is called as a non-degeneracy condition. We shall usually write things in terms of the phase function $\tau \, \Phi(x,y,t)$ of the operator $M_tf(x),$ where $\Phi(x,y,t)=x_1-y_1-A(x, y',t).$ In this case, for fixed $t,$ we write $C_t$ as 
$$C_t \, = \, \left\{(x, \tau \frac{\partial \Phi}{\partial x}, y, - \tau \frac{\partial \Phi}{\partial y}): \tau \in \mathbb{R} \setminus \{0\}, \Phi(x,y,t)=0 \mbox{ and } \psi(x,y,t) \neq 0\right\}.$$ 

In fact, the projection $\pi_l$ is a local diffeomorphism is equivalent to saying that the Jacobian of the map 
$$(\tau, y) \longrightarrow \left(\Phi(x,y,t),  \tau \frac{\partial \Phi}{\partial x} \right)$$ is non zero. The resulting Jacobian is called the Monge-Ampere determinant: 
$$J_t(x,y) \, = \, Det \left( \begin{array}{cccc}
0 & \Phi_{x_1}& \dots & \Phi_{x_n} \\
\Phi_{y_1} & \Phi_{x_1, y_1}& \dots & \Phi_{x_n,y_1} \\
\vdots & \vdots & \ddots&  \vdots\\
\Phi_{ y_n} &\Phi_{x_1, y_1} & \dots & \Phi_{x_n, y_n} \end{array} \right). $$

Since the Monge-Ampere determinant is symmetric in the $x$ and $y$ variable, we see that $\pi_l$ is a local diffeomorphism if and only if $\pi_r$ is also. Notice that, local co-ordinates can always be chosen so that we can express the full canonical relation in the form 
\begin{eqnarray} \label{4.3}
\mathcal{C} \, = \, \left\{\left(x, \xi, y, \mu,t, \tau\right): (y, \mu) \, = \, \chi_t(x, \xi), \tau= q(x,t, \xi)\right\},
\end{eqnarray}

where $q(x,t,\xi)$ is homogeneous of degree $1$ in $\xi$ and $C^{\infty}$ when $\xi \neq 0.$ To give an illustration, let us point out that for circular mean operators, we have $q(x,t, \xi)= \pm |\xi|.$ As we observed in \cite{S} that the rotational curvature condition is not sufficient  to get the local smoothing estimate for the operators $M(f)(x,t).$ Sogge showed that the following extra assumption is necessary.

Cone condition: We say that the canonical relation $\mathcal{C}$ as in \eqref{4.3} satisfies the cone condition if the cone given by the equation $\tau=q(x,t, \xi)$ has $(n-1)$ non zero principal curvatures.

It is also convenient to give a formulation which is in the spirit of the wave equation. Since $q$ is homogeneous of degree one, its Hessian with respect to the $\xi$- variable can have rank at most $n-1$. Therefore, our cone condition is: The Hessian of $q$ with respect to the $\xi$- variable has full rank, i.e., 
\begin{eqnarray} \label{4.4}
\mbox{ corank } \, q_{\xi \, \xi}^{''} \, \equiv 1.
\end{eqnarray}

If the canonical relation $\mathcal{C}$ satisfies both the non-degeneracy condition and the cone condition, we say that $\mathcal{C}$ satisfies the cinematic curvature condition (see, e.g., \cite{S}, \cite{CDS2}), since it mainly measures the way that the surfaces change with the time variable $t,$ in other words, the way in which the singularities of the operator propagate. 

In many cases, such as in the study of restriction theorems or Bochner Riesz theorems (see \cite{Z}, \cite{CS}) one wishes to prove optimal estimates for the corresponding Fourier integral operators when the mapping $\pi_{l},~\pi_{r}: C_t \rightarrow T^{\star}\mathbb{R}^n$ are allowed to be singular. To study such operators, we need the following definition due to Phong and Stein (see, \cite{Phst}).

\begin{Definition}
Let $\Sigma_t \, = \, \left\{x_0 \in C_t: \pi_l \mbox{ is } \mbox{ not } \mbox{ locally } 1-1 \right\}$ where $\pi_l, ~ \pi_r: C_t \rightarrow T^{\star}\mathbb{R}^n \setminus \{0\}$ are natural projections. We say that $C_t$ is folding of order $k-2$ if the following conditions hold: 
\begin{enumerate}
\item $\Sigma_t$ is a submanifold of $C_t$ of codimension $1$,
\item $det(d\pi_l)$ and $det(d\pi_r)$ vanish of order $k-2$ along $\Sigma_t,$
\item $T_{x_0}(\Sigma_t) \, \bigoplus \, ker (d\pi_l)_{x_0}=T_{x_0}(C_t)$,
\item $T_{x_0}(\Sigma_t) \, \bigoplus \, ker (d\pi_r)_{x_0}=T_{x_0}(C_t).$
\end{enumerate}

\end{Definition}

Last two conditions are called the transversality conditions. In fact, when $m=3,$ the above condition is equivalent to the conditions that both $\pi_l$ and $\pi_r$ are Whitney folds where they are singular (see, e.g., \cite{Hor}, \cite{SP}). $C_t$ is then called a folding canonical relation. Using this definition, Iosevich \cite{IOSE} generalized Sogge's result. He has shown that if for each $t,$ the canonical relation is folding of order $k-2$ and the cone condition is satisfied away from $\Sigma,$ then the corresponding maximal operator is $L^p$- bounded for $p >k$ in the plane. Naturally, as before, we extend this result in higher dimension.

We can now state the variable coefficient version of maximal Theorem \eqref{Theorem1.1}.
\begin{Theorem} \label{thm4.1}
Let $M_t$ be an averaging operator as in \eqref{4.1} and let $\mathcal{M}(f)(x) \, = \sup\limits_{t > 0} M_tf(x)$ be the associated maximal operator. Suppose that for each $t$ the canonical relation is folding of order $k-2$ and the cone condition \eqref{4.4} is satisfied away from $\Sigma_t.$ Then, the following inequality 
\begin{eqnarray}
\|\mathcal{M}(f)\|_{L^p(\mathbb{R}^n)} \leq \, C \, \|f\|_{L^p(\mathbb{R}^n)},
\end{eqnarray}

 holds for $p >k.$
\end{Theorem}

We now detail the dyadic decomposition of the dual space that is needed in the proof of our main result. 
\subsection{\large The dyadic decomposition} \label{subs1.1}
The proof of our main result makes use of the division of the dual (frequency) space into dyadic shells. Dyadic decomposition, whose ideas originated in the work of Littlewood and Paley, and others, will now be described in the form most suitable for us.

Let $\beta$ be a non negative radial function in $C_c^{\infty}(\mathbb{R}^n)$ supported in $\{ \frac{1}{2} \leq |\xi| \leq 2 \}$ such that 
$$\sum_{j= - \infty}^{\infty} \, \beta(2^{-j} \xi) =1 \mbox{ for } \xi \neq 0.$$
For example, we shall take,
\begin{eqnarray*}
\phi(\xi)=
\left\{\begin{array}{ll}
1, &\mbox{ if } |\xi| \leq \frac{1}{2}\\
~~0, &\mbox{ if } |\xi| \geq 1
\end{array}\right.
\end{eqnarray*}
and $$\beta(\xi) \, = \, \Phi(\frac{\xi}{2})-\Phi(\xi).$$
Then, one can easily see that  $\sum\limits_j\beta(2^{-j}\xi) \, = \, 1, ~~~\xi \, \neq \, 0$ (see \cite{JD}).

We shall use $C$ as a constant independent of $j,$ in several times without mention it. 
\section{Proof of Theorem \ref{Theorem1.1}}
Now, we proceed to prove the Theorem \ref{Theorem1.1}. Our proof will consist of three main steps. First, we shall express each operator $M_t$ as Fourier integral operator. Then we shall use the Littlewood-Paley argument and a technical lemma to reduce the problem in two dimensions. We will then use the Iosevich's approach to get our estimate for the corresponding Fourier integral operators with the help of local smoothing estimates of Mockenhaupt et. al. (see, e.g., \cite{MSS1}, \cite{MSS2}. In the process, we shall take advantage of the fact that the local smoothing argument in question is valid under small smooth perturbations.

\begin{proof} of Theorem \ref{Theorem1.1}:
We now turn to the details. First, after perhaps contracting $supp ~~(\rho),$ using a partition of unity argument, we may assume that on $supp ~~(\rho), ~~S$ can be written as   
$$S \, = \, \{x \in \mathbb{R}^n : \Phi(x)=0\},$$ where $\Phi$ is a $C^{\infty}$ function satisfying $\nabla \Phi \neq 0.$

If we use the same notation after replacing $\rho$ by $|\nabla\Phi|\rho$ and then, using the defining function $\Phi$, we can write,
\bea \label{3.2}
M_tf(x):&=&Mf(x,t)=\int\limits_{\mathbb{R}^n}\delta_0(\Phi(y)) \, f(x-ty) \, \rho(y) \, dy \nonumber \\ 
&=&\int\limits_{\mathbb{R}^n}t^{-n} \, \delta_0\left(\Phi(\frac{x-y}{t})\right) \, \rho\left(\frac{x-y}{t}\right) \, f(y) \, dy \nonumber \\
&=&\frac{1}{2\pi}\int\limits_{\mathbb{R}^n}\int\limits_{-\infty}^{\infty}t^{-n} \, e^{i\tau\Phi(\frac{x-y}{t})} \, \rho \left(\frac{x-y}{t}\right) \, f(y) \, d\tau \, dy,  
\eea

where, $\delta_0$ of course denotes  the one-dimensional Dirac-delta  function.
\noindent

Using this Fourier integral representation, we shall break up the operators dyadically. For this purpose,
let us fix $\beta \in C_c^\infty(\mathbb{R}\setminus0)$ satisfying $$\sum\limits_{-\infty}^{\infty}\beta(2^{-j}s)=1, s\neq 0.$$
We then define the dyadic operator $M_j$ by
$$M_jf(x,t)=\frac{1}{2\pi}\int\limits_{\mathbb{R}^n}\int\limits_{-\infty}^{\infty}t^{-n} \, e^{i\tau\Phi(\frac{x-y}{t})} \, \beta(2^{-j}\tau) \, \rho \left(\frac{x-y}{t}\right) \, f(y) \, d\tau \, dy.$$
Now, look at, 
\bea \label{2.2}
|M_jf(x,t)|&\leq& \frac{1}{2\pi}\int\limits_{\mathbb{R}^n}\left(\int\limits_{-\infty}^{\infty}\beta(2^{-j}\tau)d\tau\right) \, t^{-n} \, |\rho \left(\frac{x-y}{t}\right)| \, |f(y)| \, dy \nonumber \\
&\leq& \frac{1}{2\pi}2^j \int\limits_{\mathbb{R}^n}t^{-n} \, |\rho \left(\frac{x-y}{t}\right)| \, |f(y)| \, dy.
\eea

Thus, we have
\bea
\sup\limits_{t>0} \left|\sum\limits_{j\leq 0}M_jf(x,t)\right| \, &\leq& \frac{C}{2\pi}\sup\limits_{t>0}\int\limits_{\mathbb{R}^n}t^{-n} \, |\rho \left(\frac{x-y}{t}\right)| \, |f(y)| \, dy \nonumber \\
&\leq& C \,  \mathbb{M}f(x) \, = \, Hardy \, Littlewood \, maximal \, function. \nonumber
\eea

Therefore, the inequality \eqref{thm2.1} would follow from showing that when $k < p \leq \infty$, there is an $\epsilon_p > 0$ such that for $j>0,$
$$||\sup\limits_{t>0}M_jf(x,t)||_{L^p(\mathbb{R}^n)} \leq C 2^{-j\epsilon_p} \, ||f||_{L^p(\mathbb{R}^n)},$$ where $C$ is the constant independent of $j.$

Next, we claim that this is in tern would follow from 
\bea \label{2.3}
||\sup\limits_{t\in[1,2]}M_jf(x,t)||_{L^p(\mathbb{R}^n)} \leq C 2^{-j\epsilon_p}||f||_{L^p(\mathbb{R}^n)}.
\eea

\noindent 

Now to show that the inequality \eqref{2.3} is enough, 
we need to use Littlewood-paley operators $L_k$, which are defined by
$$\widehat{(L_kf)}(\xi)=\beta(2^{-k}|\xi|) \, \hat{f}(\xi).$$
Now, since $\nabla\Phi \neq 0$ on supp $\rho$, we can  see that there is an absolute constant $C_0$ such that, when $t \in [1,2]$,
\begin{eqnarray} \label{2.8}
M_jf(x,t) \, = \, M_j\left(\sum\limits_{|j-k| \leq C_0}L_kf\right)(x,t) \, + \, R_jf(x,t),
\end{eqnarray}

where, for any $N$, there is a uniform constant $C_N$ such that
$$\left|R_jf(x,t)\right| \, \leq \, C_N \, 2^{-jN} \mathbb{M}f(x),~ 1 \leq t \leq 2.$$
Thus, if \eqref{2.3} holds, then a dilation argument will give,
\begin{eqnarray*}
&&\int\sup\limits_{t>0} \, |M_jf(x,t)|^p \, dx \nonumber \\
&\leq& \sum\limits_{l=-\infty}^{\infty}\int\sup\limits_{t \in [2^l,2^{l+1}]} \, \left| M_j\left(\sum\limits_{|k+l-j|\leq C_0}L_kf\right)(x,t)\right|^p \, dx \, + \, C_N^p \, 2^{-jN_p}\int|\mathbb{M}f|^p \, dx \nonumber \\
&\leq& C^p \, C_0 \, 2^{-j \epsilon_p \, p} \, \int\sum\limits_{k=-\infty}^{\infty}|L_kf(x)|^p \, dx \, + \, C_N^p \, 2^{-jN_p}\int|\mathbb{M}f|^p \, dx \nonumber \\
&\leq& C^p \, C_0 \, 2^{-j \epsilon_p \, p}\int \left(\sum\limits_{k=-\infty}^{\infty}|L_kf|^2 \right)^{(\frac{p}{2})} \, dx \, + \, C_N^p \, 2^{-jN_p}\int|\mathbb{M}f|^p \, dx. \nonumber
\end{eqnarray*}

In the last step we have used the  fact that $p>2.$
Now using the $L^p$ boundedness of Littlewood-Paley square functions and the Hardy-Littlewood maximal theorem, we get our proof of the claim.

Now, choose a bump function $\psi \in C_c^{\infty}(\mathbb{R})$ supported in $[\frac{1}{2}, 4]$ such that $\psi(t)=1$ if $1 \leq t \leq2.$ In order to estimate \eqref{2.3}, we use the following well- known estimate (see e.g., \cite{IOSE}, Lemma 1.3),

\begin{eqnarray}
&&\sup_{t \in \mathbb{R}} \left|\psi(t) \, M_jf(x,t)\right|^p \nonumber \\
&\leq& p \, \left(\int_{-\infty}^{\infty} \left|\psi(t) \, M_jf(x,t) \right|^p \, dt \right)^{\frac{p-1}{p}} \, \left( \int_{-\infty}^{\infty} \left|\frac{\partial}{\partial t}\left[\psi(t) \, M_jf(x,t) \right]\right|^p \, dt \right)^{\frac{1}{p}}, \nonumber\\
&\leq& p \, \left(\int_{1/2}^{4} \left|M_jf(x,t) \right|^p \, dt \right)^{\frac{p-1}{p}} \, \left( \int_{1/2}^{4} \left|\frac{\partial}{\partial t}\left[M_jf(x,t) \right]\right|^p \, dt \right)^{\frac{1}{p}}, \nonumber \\
&+&  C \, p \, \int_{1/2}^{4} \left|M_jf(x,t) \right|^p \, dt \nonumber,
\end{eqnarray}

with constant $C=\|\psi'(t)\|_{L^{\infty}(\mathbb{R})}$. The first step follows by using the fundamental theorem of calculus and H$\ddot{o}$lder's inequality.
In the last step, we have used the fact that $\psi$ is supported in $[1/2, 4]$ and $\psi,~\psi'$ are uniformly bounded.

Now, integrating with respect to $x$ and by H$\ddot{o}$lder's inequality, we get,
\begin{eqnarray}
&&\|\sup_{1 \leq t \leq 2}  M_jf(x,t)\|^p_{L^p(\mathbb{R}^n)} \nonumber \\
&\leq& C \, p \, \left(\int_{\frac{1}{2}}^{4} \int_{\mathbb{R}^n} | M_jf(x,t)|^p \, dx \,  dt \right)^{\frac{p-1}{p}} \, \left( \int_{\frac{1}{2}}^{4} \int_{\mathbb{R}^n} \left|\frac{\partial}{\partial t}\left(M_jf(x,t) \right)\right|^p \,  dx \, dt \right)^{\frac{1}{p}} \nonumber \\
&& + C \, p \,  \int_{\frac{1}{2}}^{4} \int_{\mathbb{R}^n} | M_jf(x,t)|^p \, dx \,  dt. \nonumber
\end{eqnarray}

Therefore, we conclude that the inequality \eqref{2.3} would follow from 
\begin{eqnarray} \label{2.4}
&&\left( \int_{\frac{1}{2}}^{4} \int_{\mathbb{R}^n} \left|\left(\frac{\partial}{\partial t}\right)^{\alpha}M_jf(x,t)\right|^p \,  dx \, dt \right)^{\frac{1}{p}} \, \\ 
&\leq& \, C_p \, 2^{-j[\frac{1}{p} + \epsilon_p -\alpha]} \|f\|_{L^p(\mathbb{R}^n)}, ~~~ \alpha = 0, ~1. \nonumber
\end{eqnarray}

Since $\left(\frac{\partial}{\partial t}\right)M_jf(x,t)$ behaves like $2^j   \, M_jf(x,t),$ we shall only prove the estimate for $\alpha = 0.$

In proving \eqref{2.4}, we shall only use the first two co-ordinates of $y$ variables and the first two co-ordinates of $x$ variables.

With this , we set 
\begin{eqnarray}
&&\tilde{\Phi}(x,t; y) \, = \, \Phi \left( \frac{(x_1, x_2, x^{''})-(y_1, y_2, 0)}{t}\right), \nonumber \\
&& \rho(x,t; y) \, = \, \rho \left( \frac{(x_1, x_2, x^{''})-(y_1, y_2, 0)}{t}\right), \forall x^{''}=(x_3, \dots, x_n). \nonumber
\end{eqnarray}

Since, we may assume that $f$ has fixed compact support, if we define, for $g \in C_c^{\infty}(\mathbb{R}^2),$ 
\bea \label{bjexp}
B_jg(x,t) \, = \, \int_{\mathbb{R}^2} \int_{-\infty}^{\infty} e^{i \tau \tilde{\Phi}(x,t; y)} \, \beta(2^{-j} \tau) \, \rho(x,t; y) \, g(y) \, d\tau \, dy,
\eea

then it suffices to show that these operators $B_jg(x,t)$ satisfy
\begin{eqnarray}
\left( \int_{\frac{1}{2}}^{4} \int_{\mathbb{R}^n}|B_jg(x,t)|^p \,  dx \, dt \right)^{\frac{1}{p}} \, \leq \, C_p \, 2^{-j[\frac{1}{p} + \epsilon_p -\alpha]} \|g\|_{L^p(\mathbb{R}^2)}.
\end{eqnarray}

Now, we claim that, if $\epsilon_p > 0$ is small enough and if support of $\rho$ is also small, there must be a uniform constant $C$ so that we have the stronger inequality
\begin{eqnarray} \label{3.8}
&&\left( \int_{\frac{1}{2}}^{4} \int_{\mathbb{R}^2}|B_jg(x,t)|^p \,  dx_1 \, dx_2 \, dt \right)^{\frac{1}{p}} \, \\ 
&\leq& \, C \, 2^{-j[\frac{1}{p} + \epsilon_p]} \|g\|_{L^p(\mathbb{R}^2)} ~~ \forall x^{''}=(x_3, \dots, x_n). \nonumber
\end{eqnarray}

Now, we proceed to prove the inequality \eqref{3.8}. To do this, we first
notice that, for fixed $x, ~ ~\mathcal{C} \, = \, \{y \in \mathbb{R}^2: \tilde{\Phi}(x,t; y)=0\}$ denote a smooth two- dimensional curve in the plane. From our assumption on the hypersurface, we may assume that the curve $\mathcal{C}$ is a finite-type curve, which is of finite type of order $k$ at $x_0.$ Now, our proof of uniform estimates \eqref{3.8} for averaging operators associated to families of curves is based on Iosevich's approach \cite{IOSE}. We now turn into the details.

\noindent
$a)~$ We first consider the case $k=2.$

We have, 
\begin{eqnarray}
B_jg(x,t) \, &=& \, \int_{\mathbb{R}^2} \int_{-\infty}^{\infty} e^{i \tau \tilde{\Phi}(x,t; y)} \, \beta(2^{-j} \tau) \, \rho(x,t; y) \, g(y) \, d\tau \, dy \nonumber \\
&=& g \star d \mu_j, \nonumber
\end{eqnarray}

where $d \mu_j$ is the distribution of measure on the portion of the curve $\mathcal{C},$ given by 
$$d\mu_j(x-y) \, = \, \rho(x,t; y) \, \int_{-\infty}^{\infty} e^{i \tau \tilde{\Phi}(x,t; y)} \, \beta(2^{-j} \tau) \,  d\tau.$$ 

We will use the following  stationary phase result (see, \cite{CDS2}) which we shall only use for curves in $\mathbb{R}^2.$
\begin{Lemma}
Let $S$ be a smooth hypersurface in $\mathbb{R}^n$ with non-vanishing Gaussian curvature and $d \mu$ a $C_c^{\infty}$ measure on $S.$ Then
$$|\widehat{d\mu}(\xi)| \leq const \, (1+|\xi|)^{-\frac{(n-1)}{2}}.$$
Moreover, suppose that $\Gamma \subset \mathbb{R}^n \setminus \{0\}$ is the cone consisting of all $\xi$ which are normal to some point $x \in S$ belonging to a fixed relatively compact neighbourhood $N$ of supp $d\mu.$
Then,
$$\left(\frac{\partial}{\partial \xi}\right)^{\alpha} \, \widehat{d\mu}(\xi) =O\left(1+|\xi|\right)^{-N} ~ \forall N,$$ if $\xi \notin \Gamma$ and $\widehat{d\mu}(\xi) \, = \, \sum e^{-i<x_j, \xi>} \, a_j(\xi)$ if $\xi \in \Gamma,$
where the finite sum is taken over all $x_j \in N$ having $\xi$ as the normal and 
$$\left|\left(\frac{\partial}{\partial \xi}\right)^{\alpha} a_j(\xi)\right| \, \leq \, C_{\alpha} \, \left(1+|\xi|\right)^{-\frac{n-1}{2}-|\alpha|}.$$
\end{Lemma}

We now apply the lemma to each dyadic operator $B_jg(x,t)$ that is defined over a portion of our curve $\mathcal{C},$ we get an operator of the form 
$$B_jg(x,t) \, = \, \frac{1}{(2 \pi)^2} \int_{\Gamma} e^{i<x, \xi>} \, e^{it q(\xi)} \, \frac{a_j(t\xi)}{(1+t|\xi|)^{\frac{1}{2}}} \, \hat{g}(\xi) \, d\xi,$$ where $\Gamma$ is a fixed cone away from co-ordinate axes, $q(\xi)$ is homogeneous of degree one and $a_j(t\xi)$ is a symbol of order $0.$

Similarly, as before, using Littlewood-Paley square function arguments, there is an absolute constant $C_0$ such that, when $t \in [1,2],$ from  \eqref{bjexp}, we have
$$B_jg(x,t) \, = \, B_j\left(\sum\limits_{|j-k| \leq C_0}L_kg\right)(x,t) \, + \, R_jg(x,t),$$

where, for any $N$, there is a uniform constant $C_N$ such that
$$|R_jg(x,t)| \, \leq \, C_N \, 2^{-jN} \mathbb{M}g(x),~ 1 \leq t \leq 2.$$

Hence, it is enough to compute the $L^p$- norm of the following operator,
$$\tilde{B_j}g(x,t) \, = \, \sum\limits_{|j-k| \leq C_0}  \frac{1}{(2 \pi)^2} \int_{\Gamma} e^{i<x, \xi>} \, e^{it q(\xi)} \, \frac{a_j(t\xi)}{(1+t|\xi|)^{\frac{1}{2}}} \, \beta(2^{-k}|\xi|) \, \hat{g}(\xi) \, d\xi.$$

Therefore, it suffices to show that for $2 < p < \infty$ and for all $x^{''}=(x_3, \dots, x_n)$,
$$\|\tilde{B_j}g(x,t) \|_{L^p(\mathbb{R}^2\times [\frac{1}{2}, 4])}  \, \leq \, C \, 2^{-j[\frac{1}{p} + \epsilon_p^{'}]} \|g\|_{L^p(\mathbb{R}^2)}.$$

\noindent
Now, we can apply the local smoothing estimates of Mockenhaupt et. al. (see, e.g., \cite{MSS1}, \cite{MSS2}) for operators of the form
\begin{eqnarray} \label{3.9}
P_jg(x,t) \, = \,  \int e^{i<x, \xi>} \, e^{it q(\xi)} \,a(t, \xi) \, \beta(2^{-j}|\xi|) \, \hat{g}(\xi) \, d\xi,
\end{eqnarray}

 where $a(t, \xi)$ is a symbol of order $0$ in $\xi$ and the Hessian matrix of $q$ has rank $1$ everywhere.

Their results imply in particular that for $2 < p < \infty,$
\begin{eqnarray} \label{2.10}
\left( \int_{\frac{1}{2}}^{4} \int_{\mathbb{R}^2} \left|P_jg(x,t)\right|^p \, dx \, dt \right)^{\frac{1}{p}} \, \leq \, C \, 2^{j[\frac{1}{2}-\frac{1}{p} - \epsilon_p^{'}]} \|g\|_{L^p(\mathbb{R}^2)} ,
\end{eqnarray} 

for some $\epsilon_p^{'} > 0.$

Since, the operators $2^{\frac{j}{2}} \, \tilde{B_j}g(x, t)$ are of the form $P_jg(x,t),$ for suitable operators $P_j$ of the type \eqref{3.9}, 
we can apply \eqref{2.10} to obtain, 

\begin{eqnarray} \label{2.6}
&&\left( \int_{\frac{1}{2}}^{4} \int_{\mathbb{R}^2}\left|\tilde{B_j}g(x,t)\right|^p \,  dx_1 \, dx_2 \, dt \right)^{\frac{1}{p}} \, \\ 
&\leq& \, C \, 2^{-j[\frac{1}{p} + \epsilon_p]} \|g\|_{L^p(\mathbb{R}^2)} ~~ \forall x''=(x_3, \dots, x_n), \nonumber
\end{eqnarray}
if $2< p  < \infty,$ where $\epsilon_p^{'} > 0.$

However, as observed in \cite{S}, the estimate \eqref{2.10} remains valid under small, sufficiently smooth perturbations and the constant $C_p$ depends only on a finite number of derivatives of the phase function and the symbol of $P_j.$

We thus get for $2 < p < \infty$ and $\forall \,  x^{''}=(x_3, \dots, x_n)$, 
$$\|B_jg(x,t) \|_{L^p(\mathbb{R}^2 \times [\frac{1}{2}, 4])}  \, \leq \, C \, 2^{-j[\frac{1}{p} + \epsilon_p^{'}]} \|g\|_{L^p(\mathbb{R}^2)},$$
with uniform constant $C_p.$

This finishes the proof of the theorem in the case $k=2.$

\noindent
$b)~~$ Now, we consider the case $k \geq 3.$ Since our curve $\mathcal{C}$ is finite type $k,$ locally we can write $\mathcal{C}(s)=(s, \gamma(s) \, s^k +c),$ where $\gamma(s) \in C^{\infty}(\mathbb{R}), ~\gamma(0) \neq 0,$ and $c$ is a constant. The case where $k \geq 3$ can easily be reduced to the case $k=2$ by means of a dyadic decomposition in the first variable $y_1$ and re-scaling of each of the dyadic pieces in a similar way as in \cite{IOSE}. Here, we omit the details of the proof.

Hence, we finish our proof of the theorem.
\end{proof}

We give the proof of Theorem \ref{thm4.1} for the variable coefficient version of the maximal operator in the next section.

\section{Proof of Theorem \ref{thm4.1}}
Our argument will be based on a scaling argument, the proof of theorem \ref{Theorem1.1} and most importantly the local smoothing estimates of variable coefficient of Mockenhaupt, Seeger and Sogge (see, \cite{MSS1}).  As in the proof of Theorem \ref{Theorem1.1}, we localize the operator corresponding to the general family of functions by introducing a  smooth cutoff function $\psi.$ Thus, we define
$$M_t(f)(x) \, = \, \frac{1}{2\pi}\int\limits_{\mathbb{R}^n}\int\limits_{-\infty}^{\infty} \, e^{i\tau (x_1-y_1-A(x, y',t))} \, \psi(x, y, t) \, f(y) \, d\tau \, dy.$$
The only difference between the family of surfaces given by $x_1= y_1 + A(x, y',t)$ and the ones handled in Theorem \ref{Theorem1.1} is the $t$- dependence. However, using a Littlewood Paley argument as we did in the proof of Theorem \ref{Theorem1.1}, we see that it suffices to take a supremum over $t \in [1,2].$

As we did in the proof of Theorem \ref{Theorem1.1}, we first reduce the problem in two dimensions  and then it is enough to prove the $L^p$- estimate of the following operators 
$$B_j^{'}g(x,t) \, = \,  \int_{\mathbb{R}^2} \int_{-\infty}^{\infty} e^{i \tau \tilde{\Phi}(x,t; y)} \, \beta(2^{-j} \tau) \, \tilde{\psi}(x, y,t) \, g(y) \, d\tau \, dy,$$ where 
\begin{eqnarray}
&&\tilde{\Phi}(x,t; y) \, = \, x_1-y_1-  A(x, (y_2,0), t), \nonumber \\
&& \tilde{\psi}(x, y,t) \, = \, \psi (x,(y_1,y_2,0),t). \nonumber
\end{eqnarray}
Now, we are interested to get a similar estimate as \eqref{2.6} to complete the proof. To prove such estimates, our proof is exactly based on Iosevich approach \cite{IOSE} of the variable  coefficient version of maximal theorem. As before, we only need to consider the case $k=2.$
Then, we will get our result for $k \geq 3$ by means of a dyadic decomposition in the first variable and re-scaling of each of the dyadic pieces in a similar way as in \cite{IOSE}. 

The key idea of the proof is that locally our family of curves given by $y_1=x_1-A(x, (y_2,0), t),$ smoothly converges to the family of translation invariant curves with the help of  following Lemma \ref{L4.3} and an explicit proof of the lemma can be found in \cite{IOSE}.

\begin{Lemma} \label{L4.3}
Suppose that the canonical relation associated to the curve distribution $D_t$ has a $2$- sided fold of order $k-2.$ Then, for each fixed $(x^0, t')$ the curve given by the equation $y_1=x_1^{0} - A(x^0, (y_2,0), t')$ is a curve of finite type $k$ with a flat point at $y_2= \chi_{t'}(x^0),$ where $\chi_t$ is a local diffeomorphism for each fixed $(x_1, x^{''}, t),~x^{''}=(x_3, \dots, x_n)$.
\end{Lemma}

Then we can use, the local smoothing estimates to analyze the translation invariant family. And this estimate is also valid under small perturbations.
This follows from the fact that the variable coefficient estimates of Mockenhaupt et. al. \cite{MSS1} are valid under small, smooth perturbations. Such an estimate is possible, since the cinematic curvature condition is stable under smooth changes of coordinates.  In fact, we reduce our problem to the maximal operator associated to the translation invariant family, which satisfies the conclusions of Theorem \ref{thm4.1}. More precisely, we can show that if we localize our operator  corresponding to the general family of functions and take the supremum over $t$ in a sufficiently small neighbourhood of some $t',$ than the resulting operators $B_j^{'}$ maps $L^p \to L^p$ for $p>k$ with appropriate exponential decay in $j$ as required. The proof of this fact is exactly same as we did in the proof of Theorem \ref{Theorem1.1}. Here, we omit the details of the  proof.  For more details, see \cite{IOSE}. As we noted earlier, we only need to consider $t \in [1,2]$ and $(x_1, x_2)$ in some compact subset in $\mathbb{R}^2,$ for fixed $x^{''}=(x_3, \dots,x_n).$ Hence, using partitions of unity and the triangle inequality, we complete the proof.

Moreover, the proof of local smoothing estimate \cite{MSS1} shows that the estimates  depend  only on the finite number of derivatives of the phase function and the symbol of the corresponding Fourier integral operator. Hence, the estimates that are valid for the limiting operator are also valid for the sufficiently small perturbation of that operator with perhaps a larger constant. With this remark, we complete the  proof of the Theorem \ref{thm4.1}. 

\vspace{10mm}
\noindent
{\bf Acknowledgements:} 
I wish to thank the Harish-Chandra Research institute,  
the Dept. of Atomic Energy, Govt. of India, for providing excellent research facility.


\end{document}